\documentclass[11pt,a4paper]{article}
\usepackage{ifthen,latexsym,amssymb,amsmath,bbm,fixmath,xcolor}
\usepackage[shortlabels]{enumitem}

\newcommand{\C}[1]{{\protect\mathcal{#1}}}

\newcommand{\I}[1]{{\mathbbm #1}}

\renewcommand{\mid}{:}

\renewcommand{\ge}{\geqslant}
\renewcommand{\geq}{\geqslant}
\renewcommand{\le}{\leqslant}
\renewcommand{\leq}{\leqslant}

\newif\ifnotesw\noteswtrue

\newcommand{\hide}[1]{}

\newcommand{\beq}[1]{\begin{equation}\label{#1}}
\newcommand{\eeq}{\end{equation}}

\newtheorem{theorem}{Theorem}
\newtheorem{lemma}[theorem]{Lemma}
\newtheorem{proposition}[theorem]{Proposition}

\newtheorem{corollary}[theorem]{Corollary}

\newcommand{\bpf}[1][Proof.]{\smallskip\noindent{\it #1} }
\newcommand{\qed}{\nolinebreak\mbox{\hspace{5 true pt}%
  \rule[-0.85 true pt]{3.9 true pt}{8.1 true pt}}}
\newcommand{\epf}{\qed \medskip}

\newcommand{\al}{\alpha}
\newcommand{\la}{\lambda}

\newcommand{\R}{\mathbbm{R}}

\newcommand{\trace}{\mathrm{tr}}

\newcommand{\OPName}{Oleg Pikhurko
}
\newcommand{\OPAffiliation}{
Mathematics Institute and DIMAP,
University of Warwick,
Coventry CV4 7AL, UK}

\begin{document}

\newcommand{\inertia}[2]{\C I_{#1}(#2)}
\newcommand{\ratio}[2]{\C R_{#1}(#2)}

\hide{
\section*{Macros/etc}

$\inertia kG$ \verb$\inertia kG$: ACF inertia bound on $\alpha_k(G)$

$\ratio kG$ \verb$\ratio kG$: ACF ratio bound on $\alpha_k(G)$ for regular $G$

Put comma for matrix entries: $A_{i,j}$

\newpage
}

\title{A note on the Ratio and Inertia Bounds\\ for the $k$-Independence Number}
\author{Jun Gao\thanks{\OPAffiliation}\and Jie Ma\thanks{
School of Mathematical Sciences, University of Science and Technology of China, Hefei,
Anhui 230026, and Yau Mathematical Sciences Center, Tsinghua University, Beijing 100084,
China} \and \OPName\footnotemark[1]}
\maketitle

\begin{abstract}
The \emph{$k$-th power $G^k$} of a graph $G$ is the graph on the same vertex set where the edge set consists of those pairs of distinct vertices of $G$ that are at distance at most $k$ from each other.
A.~Abiad, G.~Coutinho, and M.~A.~Fiol [\emph{On the $k$-independence number of graphs}, Discrete
Mathematics 342 (2019), 2875--2885] proposed extensions of the classical ratio (for regular graphs) and inertia bounds to the independence number of $G^k$ for $k\ge 2$. 

Continuing a line of work comparing these two parameters with other known bounds, we show that the $\vartheta$-function of L.~Lov\'asz and the weighted inertia bound of A.~R.~Calderbank and P.~Frankl, when applied directly to $G^k$, perform at least as well as the ratio and inertia bounds of Abiad--Coutinho--Fiol, respectively.
In particular, $\vartheta(G^k)$ provides a polynomial-time computable upper bound on the independence number of $G^k$ that is at least as strong as the ratio bound when the latter applies (i.e.,\ when the graph $G$ is regular).
\end{abstract}

\section{Introduction}

Let $G=(V,E)$ be a graph with vertex set $[n]:=\{1,\dots,n\}$ and $k$ be a positive integer. The \emph{$k$-independence number} of~$G$, written $\al_k(G)$, is the size of a largest subset $S\subseteq V$ such that any two distinct vertices of~$S$ are at distance greater than~$k$ in~$G$. Equivalently, $\al_k(G)=\al(G^k)$ is  the standard independence number of the \emph{$k$-th power} $G^k$ of $G$, which is the graph on the same vertex set where two distinct vertices~$u$ and~$v$ are joined if they are at distance at most $k$ in $G$. The parameter $\alpha_k$ naturally arises in coding theory: in many natural ambient discrete metric spaces $(\C{X},d)$ one can build a graph $G$ on $\C{X}$ such that an error correcting/detecting code of minimum distance $k+1$ corresponds exactly to a $k$-independent set in $G$.

Abiad, Coutinho and Fiol~\cite{AbiadCoutinhoFiol19} introduced two upper bounds on $\alpha_k(G)$, which are defined as follows. For a symmetric matrix $X\in\I R^{n\times n}$, let  $\la_1(X)\ge\dots\ge\la_n(X)$ denote its eigenvalues. Let $A\in\{0,1\}^{n\times n}$ be the adjacency matrix of $G$. For $p\in\I R_k[x]$ (that is, $p$ is a real polynomial of degree at most $k$), define
\beq{eq:wpWp}
W(p):=\max_{u\in V}(p(A))_{u,u},\qquad w(p):=\min_{u\in V}(p(A))_{u,u},
\qquad \la(p):=\min_{i\in[2,n]} p(\la_i(A)).
\eeq

Define $\inertia kG$ as the minimum over all $p\in\I R_k[x]$ of 
\[
\C I^{p}_k(G):=\min\bigl\{\,|\{i: p(\la_i(A))\ge w(p)\}|,\ |\{i: p(\la_i(A))\le W(p)\}|\,\bigr\}.
\]
 Furthermore, if $G$ is regular then define $\ratio kG$ as the infimum of 
\begin{equation}
\label{eq:ACF-ratio}
\C R^{p}_k(G):=n\cdot\frac{W(p)-\la(p)}{p(\la_1(A))-\la(p)},
\end{equation} over all $p\in \I R_k[x]$ with $p(\lambda_1(A))>\lambda(p)$.

\begin{theorem}[Inertia Bound, {Abiad, Coutinho and Fiol~\cite[Theorem~3.1]{AbiadCoutinhoFiol19}}]\label{thm:ACF-inertia}
For every\linebreak[4]\ graph $G$ and every integer $k\ge 1$, it holds that $\al_k(G)\le \inertia kG$.
\end{theorem}

\begin{theorem}[Ratio Bound, {Abiad, Coutinho and Fiol~\cite[Theorem 3.2]{AbiadCoutinhoFiol19}}]\label{thm:ACF-ratio}
For every re\-gular graph $G$ and every integer $k\ge 1$, it holds that $\al_k(G)\le\ratio kG$.
\end{theorem}

For $k=1$, we get the well-known bounds on $\alpha(G)$ due to Cvetkovi\'c~\cite{Cvetkovic71} and Hoffman~\cite{Hoffman70}, see~\cite{AbiadCoutinhoFiol19}.
The functions $\inertia kG$ and $\ratio kG$ have been actively studied, including their comparisons with other known upper bounds on $\alpha(G^k)$ for various families of graphs $G$. We refer the reader to the recent survey by Abiad, Peters and Ravagnani~\cite{AbiadPetersRavagnani26}. 

Both of these bounds require choosing a polynomial $p$. As discussed in~\cite{AbiadPetersRavagnani26}, this task (for either of the bounds) can be expressed as a mixed integer linear program (MILP). The number of Boolean variables in the provided MILPs is equal to the number of distinct eigenvalues of $G$ which can potentially be as large as the number of vertices in $G$. So this MILP is not well suited for practical calculations. 

Very recently, Abiad, van Hoesel and Michaux~\cite{AbiadHoeselMichaux26} presented a polynomial-time (in $|V(G)|$) algorithm for determining the value of the Inertia Bound $\inertia kG$ for every fixed $k$. For the Ratio Bound $\ratio kG$, the  original paper of Abiad, Coutinho and Fiol~\cite{AbiadCoutinhoFiol19} provided a closed-form expression for $k=2$. An analogous closed-form expression for $k=3$ was obtained by Kavi and Newman~\cite{KaviNewman23}. 
Fiol~\cite{Fiol20} showed that if the number of closed walks in $G$ of any given length $\ell\le k$ and starting at a vertex $u$ does not depend on $u$ then $\ratio kG$ can be computed by a linear program.
A polynomial-time algorithm that computes $\ratio kG$ in general is not known. 

A well-known parameter that upper bounds $\alpha(G)$ is the \emph{$\vartheta$-function} of Lov\'asz~\cite{Lovasz79} (who used it to determine the Shannon capacity of the 5-cycle). One of its equivalent definitions (\cite[Theorem~3]{Lovasz79}) which is convenient in this note is \begin{equation}\label{eq:Theta1}
\vartheta(G):=\min\{\lambda_1(X)\mid X\in \mathcal{L}_G\},
\end{equation}
 where $\mathcal{L}_G$ is the set of $n\times n$ symmetric matrices $X$ such that
 $X_{i,j}=1$ for every $i,j\in [n]$ that are equal or non-adjacent in~$G$.

Here we show that the value of the $\vartheta$-function on~$G^k$ is at most the Ratio Bound~$\ratio kG$.

\begin{theorem}\label{th:Theta<Ratio} For every regular graph $G$ and every integer $k\ge 1$, it holds that $\vartheta(G^k)\le \ratio kG$.
\end{theorem}


As we already remarked, $\vartheta(G)\ge \alpha(G)$ for any graph $G$.
Also, $\vartheta(G^k)$ can be computed in polynomial time  in $n$ on input $(G,k)$ as a semi-definite program with one $n\times n$-matrix variable. Thus, $\vartheta(G^k)$ provides a polynomially computable upper bound on $\alpha_k(G)$ for an arbitrary graph $G$ that, by Theorem~\ref{th:Theta<Ratio}, performs at least as well as the Ratio Bound $\ratio kG$ of Abiad, Coutinho and Fiol~\cite{AbiadCoutinhoFiol19} for regular graphs.

Next, we compare the Inertia Bound $\inertia kG$ with the strengthening of  Cvetkovi{\'c}' bound~\cite{Cvetkovic71} that was introduced by Calderbank and Frankl~\cite{CalderbankFrankl92} and is defined as follows.
A \textit{weighted adjacency matrix} of $G$ is a symmetric matrix $X\in\I R^{n\times n}$ such that $X_{i,j} = 0$ whenever $i,j\in [n]$ are equal or non-adjacent in~$G$. Note that we allow arbitrary real values for adjacent indices $i$ and~$j$.
The parameter $n_{\ge 0}(G)$, which we call the \emph{Weighted Inertia Bound} here, is defined as the minimum number of non-negative eigenvalues attained by any such matrix:
\[
n_{\ge 0}(G) := \min_{X \in \mathcal{H}_G} \left|\{ i\in [n]\mid \lambda_i(X)\ge 0\}\right|,
\]
where $\mathcal{H}_G$ denotes the set of weighted adjacency matrices of $G$. It can be derived by eigenvalue interlacing that $n_{\ge 0}(G)\ge \alpha(G)$; see~\cite{CalderbankFrankl92}. Also, the classical bound by Cvetkovi{\'c}~\cite{Cvetkovic71} is obtained by allowing $X$ to be $A$ or $-A$ only.

It was verified that $n_{\ge 0}(G)=\alpha(G)$ for all graphs with at most $10$ vertices (see the discussion in~\cite[Section~6]{ElzingaGregory10}), and Godsil~\cite{Godsil04} asked if equality always holds. However, Sinkovic~\cite{Sinkovic18} showed that there are graphs (namely the 17-vertex Paley graph or this graph minus a vertex) with $n_{\ge 0}(G)$ strictly larger than $\alpha(G)$. Examples of graphs where these two parameters are strikingly far apart were constructed by Kwan and Wigderson~\cite{KwanWigderson24} and Tang, Zhang and Elphick~\cite{TangZhangElphick25}.

The computational complexity of $n_{\ge 0}(G)$ is open but sometimes the choice of weights can be guided by symmetries of the problem, as it was done in the original paper of Calderbank and Frankl~\cite{CalderbankFrankl92} for the intersection graph on $k$-sets and Huang, Klurman and Pohoata~\cite{HuangKlurmanPohoata20} for the hypercube graph.

Our second result states that $n_{\ge 0}(G^k)$ is always at most $\inertia kG$.

\begin{theorem}\label{th:Inertia}
 For every graph $G$ and every integer $k\ge 1$, it holds that
 \[
n_{\ge 0}(G^k) \le \inertia kG.
\]
\end{theorem}

\noindent\textbf{Remark.}
We would like to thank Jiang Zhou for pointing out, after the first version of this manuscript was posted on arXiv, that Theorem~\ref{th:Theta<Ratio} can also be derived from Theorem 4.1 in \cite{jiangzhou25} (by setting $M = p(A)- \lambda(p)I$ and taking $x$ to be the all-ones vector in the formula $\vartheta(\cdot)$ therein).

\section{Preliminaries and Proofs}

Let $X \in \mathbbm{R}^{n \times n}$ be a symmetric matrix. Recall that $\lambda_1(X) \geq \lambda_2(X) \geq \dots \geq \lambda_n(X)$ denote the eigenvalues of $X$.
We say that $X$ is \textit{positive semidefinite} if $\lambda_n(X)\ge 0$ (that is, all of its eigenvalues are non-negative). For symmetric matrices $X,Y\in\I R^{n\times n}$, we say that
$X$ \emph{interlaces} $Y$ if they satisfy the following chain of inequalities:
\[
    \lambda_1(X) \geq \lambda_1(Y) \geq \lambda_2(X) \geq \lambda_2(Y) \geq \dots \geq \lambda_n(X) \geq \lambda_n(Y).
\]

\hide{
We will need the following classical result by Hermann Weyl, whose proof can be found in e.g.~\cite[Theorem 4.3.1]{HornJohnson13ma}. (Note that the eigenvalues in \cite{HornJohnson13ma} are ordered non-decreasingly, which is the opposite to the ordering used here.)

\begin{lemma}[Weyl's Inequalities]\label{lem: Weyl}
Let $X$ and $Y$ be $n \times n$ real symmetric matrices. Then, for any  $i, j\in [n]$, it holds that
\begin{eqnarray*}
    \lambda_{i+j-1}(X+Y) \leq \lambda_i(X) + \lambda_j(Y),&&\mbox{if $i+j-1 \leq n$},\\
    \lambda_{i+j-n}(X+Y) \geq \lambda_i(X) + \lambda_j(Y), && \mbox{if $i+j-n \geq 1$}.
\end{eqnarray*}
\end{lemma}

Its direct consequence is the following corollary (see e.g.~\cite[Corollary 4.3.12]{HornJohnson13ma}).
}

We will need the following two results that can be directly derived from the more general inequalities of Weyl on the eigenvalues of the sum of two Hermitian matrices (see e.g.~\cite[Theorem 4.3.1]{HornJohnson13ma}). 

\begin{corollary}[{See \cite[Corollary 4.3.12]{HornJohnson13ma}}]\label{cor: psd_perturbation}
Let $X$ and $Y$ be $n \times n$ real symmetric matrices, and suppose that $Y$ is positive semidefinite. Then
\[
    \lambda_i(X) \leq \lambda_i(X + Y),\quad\mbox{for all $i = 1, \dots, n$.}
\]
\end{corollary}

\begin{corollary}[{See~\cite[Corollary 4.3.9]{HornJohnson13ma}}]\label{cor: rank1_interlace}
Let $X$ and $Y$ be $n \times n$ real symmetric matrices. Suppose that $Y$ is positive semidefinite and has rank 1. Then  $X+Y$ interlaces $X$.
\end{corollary}

The $\vartheta$-function of Lov\'asz has many other equivalent reformulations. Here we need the following minor variation where we allow diagonal entries to assume values at least $1$ (instead of being exactly $1$). Namely, we claim that
\beq{eq:Theta=Min}
\vartheta(G)=\min\bigl\{\lambda_{1}(X)\ :\ X\in\mathcal{L}^*_G\bigr\},
\eeq
 where $\mathcal{L}^*_G$ consists of all symmetric matrices $X\in \I R^{n\times n}$ such that $X_{i,j}=1$ for all non-adjacent distinct $i,j\in [n]$ and $X_{i,i}\ge 1$ for all $i\in [n]$. (Recall that the vertex set of the graph $G$ is assumed to be $[n]$.) Note that the difference between $\mathcal{L}^*_G$ and the set $\mathcal{L}_G$ from the definition of $\vartheta(G)$ in~\eqref{eq:Theta1} is that the former set allows diagonal entries to be arbitrary reals at least 1 (instead of each forced to be 1 as in the latter set). Thus $\vartheta(G)$ is at least the right-hand side of~\eqref{eq:Theta=Min}. On the other hand, take any matrix $X\in \mathcal{L}^*_G$ and let $X'$ be obtained from $X$ by changing all diagonal entries to~$1$. Then $X'\in\mathcal{L}_G$. Also,  $\lambda_1(X')\le \lambda_1(X)$ by Corollary~\ref{cor: psd_perturbation}, since $X-X'$ is a diagonal matrix with non-negative diagonal entries $X_{i,i}-1$ for $i\in [n]$. Thus~\eqref{eq:Theta=Min} holds.

\hide{
\cite[Theorem~3]{Lovasz79} showed that $\vartheta(G)$ can be equivalently defined as the maximum value of an PSD program, namely that, 
\beq{eq:Theta=Max}
\vartheta(G) =\max\Bigl\{\sum_{i,j\in V} B_{ij}\ :\ B\in\R^{V\times V}\text{ PSD},\ \trace(B)=1,\ B_{ij}=0\text{ if }ij\in E\Bigr\},
\eeq
}

Our proof of Theorem~\ref{th:Theta<Ratio} is an adaptation of the proof by Lov\'asz~\cite{Lovasz79} of the special case $k=1$ of it (namely that, for any regular graph $G$, $\vartheta(G)$ is always at most the Hoffman Bound). 

\bpf[Proof of Theorem~\ref{th:Theta<Ratio}.] Let $G$ be any regular graph on $[n]$ with eigenvalues $\lambda_1(A)\ge\dots\ge \lambda_n(A)$. Take any $p\in\I R_k[x]$ with $p(\lambda_1(A))>\lambda(p)$ and 
define 
\beq{BetaGamma}
  \beta:=\frac{n}{p(\lambda_1(A))-\lambda(p)}>0\quad\mbox{and}\quad \gamma:= W(p).
  \eeq
  (Recall that $\lambda(p)$ and $W(p)$ are defined in~\eqref{eq:wpWp}.)
Let $X:=J-\beta(p(A)-\gamma I)$, where $J$ denotes the $n\times n$ matrix with all entries equal to~$1$ and $I$ is the identity matrix.

Let us check that $X$ satisfies all assumptions in~\eqref{eq:Theta=Min}, our minor reformulation of the Lov\'asz $\vartheta$-function. Clearly, $X$ is a symmetric matrix. Take any distinct $i,j\in [n]$ which are non-adjacent in $G^k$. 
For every integer $m\in [0,k]$, the $(i,j)$-th entry of $A^m$ is 0 since it counts the number of $m$-step walks from $i$ to $j$ in $G$ while the distance between $i$ and $j$ is strictly larger than~$k$. Thus $(p(A))_{i,j}=0$ and $X_{i,j}=1$, as desired. Also, the $i$-th diagonal entry of $X$ equals $1-\beta((p(A))_{i,i}-\gamma)$, which is at least $1$ since $\beta>0$ and $\gamma=W(p)\ge (p(A))_{i,i}$. (This consideration dictated our choice of $\gamma$.)

The matrices $J$, $I$ and $A$ share the all-ones vector $\I 1$ as an eigenvector with the corresponding eigenvalue being $n$, $1$ and $\lambda_1(A)$ respectively (where we use the regularity of $G$ in the statement involving $A$). Furthermore, on the orthogonal complement of $\I 1$, the matrices $J$ and $I$ act by scalar multiplication (by $0$ and $1$ respectively). Hence the eigenvalues of $X$ are $n-\beta(p(\lambda_1(A))-\gamma)$ (with eigenvector $\I 1$) and $-\beta(p(\lambda_i(A))-\gamma)$ for $i\in [2,n]$. Thus, by $\beta$ being positive, the maximum eigenvalue $\lambda_1(X)$ of $X$ is $\max\{\, n-\beta(p(\lambda_1(A))-W(p)),\ -\beta(\lambda(p)-W(p))\,\}$. Our choice of $\beta$ makes these two terms equal, giving that $\lambda_1(X)$ equals $\C R^{p}_k(G)$, the quantity in~\eqref{eq:ACF-ratio}. Since $p\in \I R_k[x]$ was arbitrary with $p(\lambda_1(A))>\lambda(p)$, we conclude  that $\vartheta(G^k)\le \ratio kG$.\epf

Let us turn to Theorem~\ref{th:Inertia}. For a symmetric matrix $X\in \I R^{n\times n}$ and $t\in \I R$ define
    \[
    f(X,t) := |\{ i\in [n] : \lambda_i(X) \ge t \}| + | \{ j\in [n] : X_{j,j} < t \}|.
    \] 
Let $n^*_{\ge 0}(G)$ be the minimum value of $f(X,t)$ over all $X\in\mathcal{H}^*_G$ and $t\in\I R$, where  $\mathcal{H}_G^*$ denotes the set of symmetric matrices $X\in\I R^{n\times n}$ with $X_{i,j}=0$ for distinct non-adjacent $i,j\in [n]$. Note that the difference between $\mathcal{H}_G^*$ and $\mathcal{H}_G$ is that the former set allows arbitrary diagonal entries (while the latter requires that all diagonal entries are 0). 
Our next result shows that this extra freedom does not help and the parameters $n_{\ge 0}^*$ and  $n_{\ge 0}$ coincide.

\begin{proposition}\label{pr:inertia} It holds that $n_{\ge 0}^*(G) = n_{\ge 0}(G)$ for any graph $G$.
\end{proposition}


\bpf[Proof of Proposition~\ref{pr:inertia}.] 
    Suppose that $V=[n]$. 
    Then we have
    \[n^*_{\ge 0}(G)  = \min_{t,\,X\in \mathcal{H}^*_G} f(X,t) \quad \text{and} \quad n_{\ge 0}(G)  = \min_{X\in \mathcal{H}_G}  f(X,0).\]
    If we take any optimal choice of $X$ for $n_{\ge 0}(G)$ then by considering the same matrix $X$ and $t=0$ in the definition of $n^*_{\ge 0}(G)$ we conclude that $
    n^*_{\ge 0}(G) \le n_{\ge 0}(G)$.

    Let us prove the converse inequality.
    Let $t \in \mathbbm{R}$ and $X_1\in \mathcal{H}^*_G$ be such that
    \[
    n^*_{\ge 0}(G) = f(X_1,t) .
    \]
    Let $X$ be the matrix obtained from $X_1$ by setting all diagonal entries to $t$. 

    It follows from $X- tI \in \mathcal{H}_G$ and $f(X-tI,0) = f(X,t)$ that $n_{\ge 0}(G)  \le f(X,t)$. 
    So it is enough to show $ f(X,t) \le n^*_{\ge 0}(G)$.   
    Let $X_1, X_2, \dots, X_n, X_{n+1}$ be a sequence of matrices, where for each $i \in [n]$, $X_{i+1}$ is obtained from $X_i$ by setting the $i$-th diagonal entry $(X_i)_{i,i}$ to $t$. Clearly $X= X_{n+1}$. Fix any $r\in[n]$. 
    
    If $(X_r)_{r,r} \ge t$ then, by Corollary~\ref{cor: rank1_interlace}, we have 
    \[
    |\{i: \lambda_i(X_{r+1})\ge t\}| \le |\{i: \lambda_i(X_{r})\ge t\}|
    \] 
    while clearly $|\{j:(X_{r+1})_{j,j}<t\}| = |\{j:(X_r)_{j,j}<t\}|$, which implies that $f(X_{r+1},t) \le f(X_{r},t)$.
    
    If $(X_r)_{r,r} < t$,  then we have $|\{j:(X_{r})_{j,j}<t\}| = |\{j:(X_{r+1})_{j,j}<t\}|+1$.
    Let $j$ be the maximum integer in $[n]$ such that $\lambda_j(X_r) \ge t$ (and set $j=0$ if no such integer exists), so that $|\{i: \lambda_i(X_{r})\ge t\}|=j$.
    By Corollary~\ref{cor: rank1_interlace}, $\lambda_{i+1}(X_{r+1}) \le \lambda_{i}(X_r)$ for every $i\in[n-1]$. In particular, if $j\le n-2$ then $\lambda_{j+2}(X_{r+1}) \le \lambda_{j+1}(X_r) <t$, so all eigenvalues of $X_{r+1}$ of index at least $j+2$ are smaller than $t$ and
    \[
    |\{i: \lambda_i(X_{r+1})\ge t\}| \le j+1 = |\{i: \lambda_i(X_{r})\ge t\}| +1.
    \]
    If $j\ge n-1$ then the above bound is trivial.
    So we also have $f(X_{r+1},t) \le f(X_{r},t)$.
    
    Thus for any $r\in[n]$ we have
    \[ f(X_{r+1},t) \le f(X_{r},t),\]
    which implies that $ f(X,t)=f(X_{n+1},t) \le f(X_1,t) = n^*_{\ge 0}(G)$.
\epf

\bpf[Proof of Theorem~\ref{th:Inertia}.] Take any $p\in\I R_k[x]$. 
Then both $p(A)$ and $-p(A)$ are  in $\mathcal{H}_{G^k}^*$ since $p(A)_{i,j}=0$ for any $i,j\in [n]$ at distance more than $k$ in $G$. Observe that 
\begin{eqnarray*}
f(p(A),w(p))&=&|\{i: p(\la_i(A))\ge w(p)\}|,\\
f(-p(A),-W(p))&=&|\{i: p(\la_i(A))\le W(p)\}|.
\end{eqnarray*} 
(Indeed, $w(p)=\min_u(p(A))_{u,u}$ and $W(p)=\max_u(p(A))_{u,u}$, so the second summand in the definition of $f$ vanishes in each case.) Thus 
\[
n_{\ge 0}^*(G^k)\le \min\left\{f(p(A),w(p)),\,f(-p(A),-W(p))\right\}=\C I^p_k(G).
\] 
Since $p\in\I R_k[x]$ was arbitrary, we conclude that  $n_{\ge 0}^*(G^k)\le \inertia kG$. By Proposition~\ref{pr:inertia}, $n_{\ge 0}(G^k)=n_{\ge 0}^*(G^k)\le \inertia kG$, as desired.\epf

\section*{Acknowledgements}

Jun Gao and Oleg Pikhurko were supported by ERC Advanced Grant 101020255. Jie Ma was supported by the National Key Research and Development Program of China 2023YFA1010201 and the National Natural Science Foundation of China grant 12125106. Jie Ma would like to thank the Mathematics Institute and DIMAP at the University of Warwick for their hospitality during his visit.

AI (Claude Opus 4.8) was used to proofread the final version of this paper.

\bibliography{paper}

@article{TangZhangElphick25,
  title={Inertia, independence and expanders},
  author={Q. Tang and S. Zhang and C. Elphick},
  journal={Bulletin of the London Mathematical Society},
  volume={57},
  pages={4076--4095},
  year={2025}
}

@article{Fiol20,
  title={A new class of polynomials from the spectrum of a graph, and its application to bound the {$k$}-independence number},
  author={M. A. Fiol},
  journal={Linear Algebra and its Applications},
  volume={605},
  pages={1--20},
  year={2020}
}

@article{HuangKlurmanPohoata20,
  author  = {Huang, H. and Klurman, O. and Pohoata, C.},
  title   = {On subsets of the hypercube with prescribed {H}amming distances},
  journal = {Journal of Combinatorial Theory, Series A},
  volume  = {171},
  year    = {2020},
  pages   = {105156, 21 pp.},
  doi     = {10.1016/j.jcta.2019.105156},
}

@article{Sinkovic18,
  author  = {Sinkovic, J.},
  title   = {A graph for which the inertia bound is not tight},
  journal = {Journal of Algebraic Combinatorics},
  volume  = {47},
  year    = {2018},
  pages   = {39--50},
}

@article{ElzingaGregory10,
  title={Weighted matrix eigenvalue bounds on the independence number of a graph},
  author={R. J. Elzinga and D. A. Gregory},
  journal={Electronic Journal of Linear Algebra},
  volume={20},
  pages={468--489},
  year={2010}
}

@unpublished{Godsil04,
  author = {Godsil, C.},
  title  = {Interesting Graphs and their Colourings},
  year   = {2004},
  note   = {Unpublished lecture notes},
}

@article{KwanWigderson24,
  author  = {Kwan, M. and Wigderson, Y.},
  title   = {The inertia bound is far from tight},
  journal = {Bulletin of the London Mathematical Society},
  volume  = {56},
  year    = {2024},
  pages   = {3196--3208}
  }

@article{AbiadCoutinhoFiol19,
  author =        {A. Abiad and G. Coutinho and M. A. Fiol},
  journal =       {Discrete Mathematics},
  pages =         {2875--2885},
  title =         {On the {$k$}-independence number of graphs},
  volume =        {342},
  year =          {2019},
}

@article{Cvetkovic71,
  author =        {D. M. Cvetkovi{\'c}},
  journal =       {Publikacije Elektrotehni\v{c}kog Fakulteta.
                   Univerzitet u Beogradu. Serija: Matematika i Fizika},
  pages =         {1--50},
  title =         {Graphs and their spectra},
  volume =        {354--356},
  year =          {1971},
}

@incollection{Hoffman70,
  address =       {New York},
  author =        {A. J. Hoffman},
  booktitle =     {Graph Theory and its Applications},
  editor =        {B. Harris},
  pages =         {79--91},
  publisher =     {Academic Press},
  title =         {On eigenvalues and colorings of graphs},
  year =          {1970},
}

@article{AbiadPetersRavagnani26,
  author =        {A. Abiad and L. Peters and A. Ravagnani},
  journal =       {Canadian Journal of Mathematics},
  note =          {Published online;
                   \texttt{doi:10.4153/S0008414X25101600}},
  title =         {The {E}igenvalue {M}ethod in coding theory},
  year =          {2026},
}

@unpublished{AbiadHoeselMichaux26,
  author =        {A. Abiad and S. van Hoesel and V. Michaux},
  note =          {arXiv:2605.00524},
  title =         {Optimization and complexity of inertia-type bounds on
                   the independence and chromatic numbers of graph
                   powers},
  year =          {2026},
}

@article{KaviNewman23,
  author =        {L. C. Kavi and M. Newman},
  journal =       {Discrete Mathematics},
  pages =         {113471},
  title =         {The optimal bound on the 3-independence number
                   obtainable from a polynomial-type method},
  volume =        {346},
  year =          {2023},
}

@article{Lovasz79,
  author =        {L. Lov{\'a}sz},
  journal =       {IEEE Transactions on Information Theory},
  pages =         {1--7},
  title =         {On the {S}hannon capacity of a graph},
  volume =        {25},
  year =          {1979},
}

@article{CalderbankFrankl92,
  author =        {A. R. Calderbank and P. Frankl},
  journal =       {Combinatorics, Probability and Computing},
  pages =         {115--122},
  title =         {Improved upper bounds concerning the
                   {E}rd{\H{o}}s--{K}o--{R}ado theorem},
  volume =        {1},
  year =          {1992},
}

@book{HornJohnson13ma,
  author =        {Horn, R. A. and Johnson, C. R.},
  edition =       {Second},
  pages =         {xviii+643},
  publisher =     {Cambridge Univ.\ Press},
  title =         {Matrix analysis},
  year =          {2013},
}

@article {jiangzhou25,
    author = {Zhou, Jiang},
     title = {Unified bounds for the independence number of graphs},
  journal = {Canadian Journal of Mathematics.},
    volume = {77},
      year = {2025},
    number = {1},
     pages = {97--117},
       DOI = {10.4153/S0008414X23000822},
       URL = {https://doi.org/10.4153/S0008414X23000822},
}
\end{document}